\documentclass{article}
\usepackage{amsmath}
\usepackage{amsfonts}
\usepackage{amssymb}
\usepackage{latexsym}
\usepackage{epsfig}
\usepackage{colordvi}
\usepackage{graphicx}
\usepackage{color}

\newcommand{\BBH}{\mbox{$\mathbb{H}$}}
\newcommand{\BBR}{\mbox{$\mathbb{R}$}}

\newcommand{\BBV}{\mbox{$\mathbb{V}$}}

\newcommand{\bD}{\mbox{\boldmath{$D$}}}

\newcommand{\fb}{\mbox{\boldmath{$f$}}}

\newcommand{\bI}{\mbox{\boldmath{$I$}}}

\newcommand{\bn}{\mbox{\boldmath{$n$}}}

\newcommand{\bq}{\mbox{\boldmath{$q$}}}

\newcommand{\bu}{\mbox{\boldmath{$u$}}}

\newcommand{\sbv}{\mbox{\boldmath{\small $v$}}}
\newcommand{\bv}{\mbox{\boldmath{$v$}}}
\newcommand{\bw}{\mbox{\boldmath{$w$}}}

\newcommand{\bx}{\mbox{\boldmath{$x$}}}

\newcommand{\bsigma}{\mbox{\boldmath{$\sigma$}}}

\newcommand{\sbtau}{\mbox{\boldmath{\footnotesize $\tau$}}}

\newcommand{\bpsi}{\mbox{\boldmath{$\psi$}}}

\newcommand{\bzero}{\mbox{$\bf 0$}}

\newcommand{\dis}{\displaystyle}

\makeatletter\@addtoreset{equation}{section}\makeatother

\newtheorem{theorem}{Theorem}[section]

\newtheorem{corollary}{Corollary}[section]

\newtheorem{definition}{Definition}[section]

\newcommand{\eps}{\varepsilon}

\DeclareMathOperator{\Div}{div}

\newcommand{\QED}{\hfill\ensuremath{\square}}

\begin{document}

\begin{center}
{\Large {\sc Analysis of the Brinkman-Forchheimer equations with slip boundary conditions}}
\end{center}

\begin{center}

J. K. Djoko\footnote{Is the corresponding author. jules.djokokamdem@up.ac.za}~\quad\&\quad~P. A. Razafimandimby\\
Department of Mathematics and Applied Mathematics,\\
University of Pretoria, \\
Private bag X20, Hatfied 0028, Pretoria, South Africa~

\end{center}
\makebox[1cm]{ }\vspace{3ex}\\

{\bf Abstract}\\
In this work, we study the Brinkman-Forchheimer equations driven under
slip boundary conditions of friction type. We prove the existence
and uniqueness of weak solutions by means of regularization
combined with the Faedo-Galerkin approach. Next we discuss the
continuity of the solution with respect to Brinkman's and
Forchheimer's coefficients. Finally, we show that the weak solution
of the corresponding stationary problem is stable.



\section{Introduction}

We consider the Brinkman-Forchheimer equations for unsteady flows of
incompressible fluids, i.e.
\begin{eqnarray}
\label{eq:1-1a}
\dis\frac{\partial \bu}{\partial t}-\nu\Delta \bu+a\bu +b|\bu|^\alpha \bu+\nabla p&=&\fb\,\,\text{ in } Q=\Omega \times (0,T)\,,\\
\label{eq:1-1b} \Div \bu&=&0\,\,\text{ in } Q\,,
\end{eqnarray}
where $\Omega$ is the flow region, a bounded domain in $\BBR^3$. The
motion of our incompressible fluid is described by the velocity
$\bu(\bx,t)$ and pressure $p(\bx,t)$. In (\ref{eq:1-1a}) and
(\ref{eq:1-1b}), $\fb$ is the external body force per unit volume
depending on $\bx$ and $t$, while the positive parameters $\nu,a,b$
are respectively the Brinkman coefficient, the Darcy coefficient and
Forchheimer coefficient, and $\alpha\in[1,2]$ is a given number.
Equations (\ref{eq:1-1a}) and (\ref{eq:1-1b}) are supplemented by
boundary and initial conditions. As far as the initial condition
goes, we assume that
\begin{equation}
\bu(\cdot,0)=\bu_0~~\mbox{on}~\overline{\Omega}~, \label{eq:1-2}
\end{equation}
where $\bu_0$ is a given function, that will be made precise later,
and $\overline{\Omega}$ is the closure of $\Omega$. Next in order to
describe the motion of the fluid at the boundary,  we assume that
the boundary of $\Omega$, say, $\partial \Omega$ is made of two
components $S$ (say the outer wall) and $\Gamma$ (the inner wall),
and it is required that $\overline{\partial\Omega} =
\overline{S\cup\Gamma}$, with $S\cap\Gamma=\emptyset$. We assume the
homogeneous Dirichlet condition on $\Gamma$, that is
\begin{equation}
\label{eq:1-3} \bu  = 0\quad \mbox{on}~\Gamma\times(0,T)~.
\end{equation}
We have chosen to work with homogeneous condition on the velocity in
order to avoid the technical arguments linked to the Hopf lemma (see
\cite{Girault-Raviart}, Chapter 4, Lemma 2.3)\,. On $S$, we first
assume the impermeability condition
\begin{equation}
u_N=\bu\cdot\bn=0\,\,\,\text{on}\,\, S\times(0,T)\,, \label{eq:1-4}
\end{equation}
where $\bn$ is the outward unit normal on the boundary $\partial
\Omega$, and $u_N$ is the normal component of the velocity, while
$\bu_{\sbtau} =\bu - u_N\bn$ is its tangential component. In
addition to  (\ref{eq:1-4}) we also impose on $S$, a threshold slip
condition \cite{Leroux,Fujita94}, which is the main ingredient of
this work. The threshold slip condition can be formulated with the
knowledge of a positive function $g:S\longrightarrow (0,\infty)$
which is called the barrier of threshold function and the use of
sub-differential to link quantities of interest. It is written as
\begin{equation}
-(\bsigma\bn)_{\sbtau}\in~g\partial
|\bu_{\sbtau}|\quad\mbox{on}~S\times(0,T) \,, \label{eq:1-5}
\end{equation}
where $(\bsigma\bn)_{\sbtau}$ is the tangential component of the
Cauchy tensor $\bsigma$ given by $\bsigma =-p\bI + 2\nu\bD(\bu)$
with $\bD(\bu)=\frac{1}{2}[\nabla\bu + (\nabla\bu)^T]$, and
$\partial|\cdot|$ is the sub-differential of the real valued
function $|\cdot|$, with $|\bw|^2=\bw\cdot\bw$\,. We recall that if
$X$ is a Hilbert space with $x_0\in X$, then
\begin{equation}
y\in \partial \Psi(x_0)~\mbox{if and only if}~\Psi(x)-\Psi(x_0)\geq
y\cdot(x-x_0)\quad\forall x\in X\,. \label{eq:definition}
\end{equation}
Without using the sub-differential, the threshold condition
(\ref{eq:1-5}) can be written as \cite{DUVAUT}
\begin{equation}
\label{eq:slip} \left.
\begin{aligned}
|(\bsigma\bn)_{\sbtau}| &\leq g,   \\
|(\bsigma\bn)_{\sbtau}| &< g  \Rightarrow  \bu_{\sbtau} = \bzero,
\\
|(\bsigma\bn)_{\sbtau}| &= g  \Rightarrow \bu_{\sbtau} \neq
\bzero~,~
-(\bsigma\bn)_{\sbtau}=g\frac{\bu_{\sbtau}}{|\bu_{\sbtau}|}\\
\end{aligned} \ \right\} \ \mbox{on}\quad S\times (0,T)~.
\end{equation}
One observes that different boundary conditions describe different
physical phenomena. In \cite{Rajagopal}, the equations of Brinkman
corresponding to (\ref{eq:1-1a}) with $b=0$  have been derived using
mixtures theory, in fact a class of approximate models for flows of
incompressible fluids passing porous solids have been described.
{\it Forchheimer} \cite{Forchheimer} studied flow experiments in
sandpacks and came to the conclusion that for small Reynolds numbers
(Re $\approx$ 0.2), the diffusion law of Darcy is not significant.
Furthermore, he found the relationship between the pressure gradient
and the velocity obtained using the law of Darcy to be nonlinear. In
fact for a wide range of physical experiments, he found that the
nonlinear term should be quadratic. Inertial effects in the porous
medium at relatively small Reynolds numbers are the cause of the
nonlinear excess pressure drop observed by {\it Forchheimer} and
others. The slip boundary conditions of friction type (\ref{eq:1-5})
can be justified by the fact that frictional effects of the fluid at
the pores of the solid can be very important. In fact on the role of
the boundary conditions for such problems, {\it Brinkman}
\cite{Brinkman} mentioned that ``The flow through this porous media
is described by a modification of Darcy's equation. Such
modification was necessary to obtain consistent boundary
conditions''. While there continues to be some debate over the
functionality of the Brinkman-Forchheimer model \cite{Nield},
nonlinearity has been verified experimentally \cite{Kladias-Prasad},
and some theoretical results have been obtained in
\cite{Celebi,Payne,AmesStraughan,AmesPayne,Payne2}. The
Brinkman-Forchheimer equation continues to be used for predicting
the velocity of the flow in the porous region, while the energy
equation for the porous region accounts for the effect of thermal
dispersion \cite{Kuznetsov-Xiong}. Since we are well aware that for
such flow, there are important features at the boundary, it is
therefore important to model Brinkman-Forchheimer flow accurately
taking into account the motion at the boundary. This is the driving
force behind our work.\newline Even though flows under boundary
conditions of friction type have been considered in various
publications
(\cite{Leroux,Fujita94,LerouxTani,FujKawSas95,FujKawKaw95,FujKaw98,Fujita01,
Fujita02,Fujita02b} among others), and Brinkman-Forchheimer
equations (\ref{eq:1-1a}), (\ref{eq:1-1b}) with non slip boundary
conditions has been examined in
\cite{Celebi,Payne,AmesStraughan,AmesPayne,Payne2}, the combination
of (\ref{eq:1-1a}), (\ref{eq:1-1b}) and (\ref{eq:1-5}) has not been
presented in the literature, and it is the object of this work. The
novelty of the problem, from the mathematical point of view, derives
from the following features; the highly coupled and nonlinear nature
of the problem, the incompressibility constraint and related
pressure, and the leak boundary conditions (\ref{eq:1-4}) and
(\ref{eq:1-5}).\newline Not surprisingly, flows problems involving
boundary conditions of friction type offer significant theoretical
and computational challenges. With regard to theoretical studies,
the work by Hiroshi Fujita and co-authors
\cite{Fujita94,FujKawSas95,FujKawKaw95,FujKaw98,Fujita01,
Fujita02,Fujita02b}, represent some early, contributions. These
authors established existence, and uniqueness of solutions, for the
equations corresponding to Stokes equations by means of semi-group
approach, regularity of solutions are also examined. An interesting
and related work is that by Christiaan Leroux and co-author
\cite{Leroux,LerouxTani} on Stokes and Navier Stokes equations under
more general ``friction type boundary conditions''. As far as
computational studies for flows driven by ``friction type boundary
conditions'' are concerned, it should be mentioned that even though
the literature is very rich in problems formulated in terms of
variational inequalities
\cite{GlowinskiLionsTremoliere,Glowinski,WeiminReddy,SHS06}, not
much have been done for the specific case involving mixed coupled
problems \cite{WeiminReddy,Yuan-Kaito,GerLel09,Sasamoto92,SasKaw93},
and we would like to explore that research direction.\newline
Problem (\ref{eq:1-1a})-(\ref{eq:1-5}) is a coupled nonlinear system
of equations with a non-differentiable expression (at zero) on $S$
due to the sub-differential term $\partial |\bu_{\sbtau}|$. We
propose to solve the resulting system of partial differential
equations using the regularization approach \cite{DUVAUT,LIONS},
which consists of replacing the original problem by a sequence of
``better behaved'' approximate problems indexed by a small positive
parameter $\eps$. We then solve the regularized problems by the
Faedo-Galerkin method, and finally, the solution of the original
problem is obtained by passage to the limit as $\eps$ goes to zero.
The difficulty in the algorithm described is to obtain the pressure.
Indeed, as the problem in its weak form is formulated as a
variational inequality with only one unknown, the pressure will not
be obtained in the usual way (for the classical Navier-Stokes
equations see e.g \cite{Raviart-Thomas}, [Theorem 2.5-1, page 54]).
But, instead we first construct a regularized pressure by using the
classical approach and then pass to the limit as $\eps$ goes to
zero, after showing that the regularized pressures are bounded in
some appropriate function space. After constructing weak solutions
of the problem, we analyze some qualitative properties of the
solution, namely; the continuous dependence of the solution with
respect to the {\it Brinkman} and {\it Forchheimer} coefficients,
and the stability of the stationary solution. The results presented,
extend in some sense those obtained in \cite{Payne,Payne2} to a
family of variational inequalities with non-differentiable
functionals.\newline The remaining part if this work is organized as
follows. In section 2, we document the variational formulation
associated to the problem and prove its well-posedness. Section 3 is
devoted to the stability of the solutions with respect to some data
of the problem. The stability of the stationary solutions is
analyzed in Section 4.

\section{Analysis of the problem: Solvability}

We introduce some preliminaries and notation for the mathematical
setting of the problem. We write down a variational formulation of
problem (\ref{eq:1-1a})-(\ref{eq:1-5}). Next we derive some {\it a
priori} estimates of its solution and obtain the existence of
solutions by means of Faedo-Galerkin.

\subsection{Preliminaries/Notation}

In what follows, for $1\leq p\leq \infty$, $L^p(\Omega)$, and
$L^p(\partial \Omega)$ are the usual Lebesgue spaces, with norms
denoted by $\|\cdot\|_{L^p}$ and $\|\cdot\|_{L^p(\partial \Omega)}$
respectively. (of course when $p=2$, we will denoted the norm in
$L^2(\Omega)$ by $\|\cdot\|$). We shall use the following notation;
for the sake of simplicity, one defines them in three dimensions.
Let $k=(k_1,k_2,k_3)$ denote a triple of non-negative intergers, set
$|k|=k_1+k_2+k_3$ and define he partial derivative $\partial^k$ by
\[
\partial^kv=\dis\frac{\partial^{|k|}v}{\partial x^{k_1}\partial y^{k_2}\partial
z^{k_3}}~.
\]
Then, for non-negative integer $m$, we recall the classical Sobolev
space
\[
H^m(\Omega)=\{v\in L^2(\Omega)~;~~\partial^kv\in
L^2(\Omega)~~\forall~|k|\leq m\}
\]
equipped with the seminorm
\[
|v|_{H^m(\Omega)}=\left[\sum_{|k|=m}\int_\Omega
|\partial^kv|^2\mbox{dx}\right]^{1/2}
\]
and norm
\[
||v||_{H^m(\Omega)}=\left[\sum_{0\leq k\leq m
}\int_\Omega|\partial^kv|^{2}\mbox{dx}\right]^{1/2}~.
\]
For $p=1,2,3, \cdots$, the inner products in the spaces
$L^2(\Omega)^p$, $L^2(\partial \Omega)^p$ and $H^1(\Omega)^p$ are
denoted by $(\cdot,\cdot)$, $(\cdot,\cdot)_{\partial \Omega}$ and
$(\cdot,\cdot)_1$, respectively. The product spaces are denoted by
bold letters: $\mathbf{H}^1(\Omega)=H^1(\Omega)^3$,
$\mathbf{L}^2(\Omega)=L^2(\Omega)^3$,
$\mathbf{L}^{\alpha+2}(\Omega)=L^{\alpha+2}(\Omega)^3$, etc.\newline
Here, and in what follows, the boundary values are to be understood
in the sense of traces. We omit the trace operators where the
meaning is direct; otherwise we denote the traces by $\bv|_\Gamma$,
$\bv|_S$, etc. Also, all the derivatives should be understood in the
sense of distribution.\newline We also recall from
\cite{Girault-Raviart} (Chap. I, Thm 1.1) for instance the following
Poincar\'{e}-Friedrichs inequality:
\begin{equation}
\text{for all}\quad \bv\in \mathbf{H}^1(\Omega)
\cap\{v_n|_S=0~,~\bv|_\Gamma=0\}\,,\quad \|\bv\|\leq
C\|\nabla \bv\|\,\,, \label{eq:Poincare}
\end{equation}
which yields the equivalence of the norms $\|\cdot\|_1$ and
$|\cdot|_1$ on
$\mathbf{H}^1(\Omega)\cap\{v_n|_S=0~,~\bv|_\Gamma=0\}$.\newline For
any separable Banach space $E$ equipped with the norm $\|\cdot\|_E$,
we denote by $C^0(0,T;E)$ the space of continuous functions from
$[0,T]$ with values in $E$ and by $D'(0,T;E)$ the space of
distributions with values in $E$. $L^p(0,T;E)$ is a Banach space
consisting of (classes of) functions $t\longmapsto f(t)$ measurable
from $[0,T]\longmapsto E$ (for the measure $dt$) such that
\begin{eqnarray*}
&&\|f\|_{L^p(0,T;E)}=\left[\int^T_0\|f(t)\|^p_Edt\right]^{1/p}<\infty~\text{for}~p\neq \infty\\
&&\|f\|_{L^\infty(0,T;E)}=\dis\text{ess}_{0<t<T}\sup\|f(t)\|_E<\infty~.
\end{eqnarray*}
In what follows, $\phi(t)$ stands for the function $\bx\in
\Omega\mapsto \phi(\bx,t)$. \newline We assume that the data
$(\fb,g)$ belong to $L^2(0,T;\mathbf{L}^2(\Omega))\times
L^\infty(S)^2$, and that the datum $\bu_0$ belongs to
$\mathbf{H}^1(\Omega)\cap \mathbf{L}^{\alpha+2}(\Omega)$, and
satisfies the incompressibility condition
\begin{equation}
\Div\bu_0=0\quad\text{in}\quad\Omega\,. \label{eq:compatibility}
\end{equation}
This last condition is not necessary for all the results that
follow but, since it is not restrictive, we shall assume it from
now on.

\subsection{Variational formulation}

In order to write a variational form associated with
(\ref{eq:1-1a})-(\ref{eq:1-5}), we retain (\ref{eq:1-2}) and we
weaken the equations (\ref{eq:1-1a}), (\ref{eq:1-1b}) and
constraints (\ref{eq:1-3}), (\ref{eq:1-4}) using the Green's
formula, while (\ref{eq:1-5}) is re-interpreted with the help of
(\ref{eq:definition}). It follows from the nonlinear term in
(\ref{eq:1-1a}) that $\bu(t)$ and the test function $\bv$ should
belong to $\mathbf{L}^{\alpha+2}(\Omega)$. Then $\bu'(t)$ and
$|\bu(t)|^\alpha\bu(t)$ must belong to the conjugate of
$\mathbf{L}^{\alpha+2}(\Omega)$, which is
$\mathbf{L}^{\frac{\alpha+2}{\alpha+1}}(\Omega)$\,. We then
introduce the following spaces
\begin{eqnarray*}
{\mathcal N}&=&\mathbf{H}^1(\Omega)\cap\{\bv|_\Gamma=0\,\,,\,v_n|_S=0\}~,\nonumber\\
M&=&L^2_0(\Omega)=\{q\in L^2(\Omega)~,~(q,1)=0\}~.\nonumber
\end{eqnarray*}
We then adopt the following definition of weak solutions of
(\ref{eq:1-1a})--(\ref{eq:1-5})
\begin{definition}
Given $(\fb,g)$ in $L^2(0,T;\mathbf{L}^2(\Omega))\times
L^\infty(S)^2$, and $\bu_0 \in \mathbf{H}^1(\Omega)\cap
\mathbf{L}^{\alpha+2}(\Omega)$, satisfying (\ref{eq:compatibility}).
We say that $(\bu,p)$ is a weak solution of
(\ref{eq:1-1a})--(\ref{eq:1-5}) if and only if;\newline $\bu\in
L^\infty\left(0,T;({\mathcal N}\cap
\mathbf{L}^{\alpha+2}(\Omega)\right)$, $p\in L^2(0,T;M)$~, and
$\bu'\in L^2\left(0,T;({\mathcal N}\cap
\mathbf{L}^{\alpha+2}(\Omega))'\right)$, and\newline for almost all
$t$ and all $~q \in L^2(\Omega), \bv\in {\mathcal N}\cap
\mathbf{L}^{\alpha+2}(\Omega)$
\begin{eqnarray}
&&(\bu'(t),\bv-\bu(t))+\gamma
(\nabla\bu(t),\nabla(\bv-\bu(t))\,)+a\,(\bu(t),\bv-\bu(t))\nonumber\\
&&+b\,(|\bu(t)|^\alpha \bu(t),\bv-\bu(t))
-(\Div(\bv-\bu(t)),p(t))+\nonumber\\
&&+J(\bv)-J(\bu(t)) \geq (\fb(t),\bv-\bu(t))\,,\label{eq:2-1}\\
&&(\Div\bu(t),q)=0\,\,\label{eq:2-2}\\
&&\bu(0)=\bu_0~,\label{eq:2-3}
\end{eqnarray}
where, $J(\bv)=(g(\bx),|\bv_{\sbtau}(\bx)|)_{L^2(S)}$\,.
\label{weaksolution}
\end{definition}
Following \cite{DUVAUT}, it can be shown that any solution of
(\ref{eq:1-1a})-(\ref{eq:1-5}) is a solution of
(\ref{eq:2-1})-(\ref{eq:2-3}) in the sense of distributions. The
converse property holds for any solution of the problem
(\ref{eq:1-1a})-(\ref{eq:1-5}) that enjoys the regularity mentioned in 
Definition \ref{weaksolution}, in a sense to be made precise later
on. The kernel of the bilinear and continuous form $
L^2(\Omega)\times {\mathcal N}\cap
\mathbf{L}^{\alpha+2}(\Omega)\ni(q,\bv) \longmapsto
(q,\Div\bv)\in\BBR$~ is $\BBV=\{\bv\in {\mathcal N}\cap
\mathbf{L}^{\alpha+2}(\Omega)\,\,,\Div\bv=0\quad\text{in}\quad\Omega\}$\,.
With the space $\BBV$ in mind, it is then easy to see that the
function $\bu(t)$ given in (\ref{eq:2-1})--(\ref{eq:2-3}) is a
solution of the simpler variational problem: Find $\bu\in
L^\infty(0,T;\BBV)~,~\bu'\in L^2(0,T;\BBV')$ satisfying
(\ref{eq:2-3}) such that for almost all $t$ and all~$\bv\in \BBV$
\begin{eqnarray}
&&(\bu'(t),\bv-\bu(t))+\gamma
(\nabla\bu(t),\nabla(\bv-\bu(t))\,)+a\,(\bu(t),\bv-\bu(t))~\label{eq:2-4}\\
&&+ b\,(|\bu(t)|^\alpha
\bu(t),\bv-\bu(t))+J(\bv)-J(\bu(t))\geq (\fb(t),\bv-\bu(t))\,.\nonumber
\end{eqnarray}
Next, we establish the solvability of the variational problem
(\ref{eq:2-4}) by means of regularization combined with Galerkin's
method. We then construct a pressure $p$ in $L^2(0,T;L^2_0(\Omega))$
such that the couple $(\bu,p)$ enjoys the regularity announced in
definition \ref{weaksolution}, and satisfies
(\ref{eq:2-1})--(\ref{eq:2-3})\,.

\subsection{Existence of a solution}

In this paragraph we discuss the solvability of (\ref{eq:2-4}) by
regularization, and passage to the limit. Thus it is obtained in
several steps, that we describe below.

Step 1: {\it Regularized problem}.\newline We first recall that one
of the difficulties of solving (\ref{eq:2-4}) is the fact that the
functional $\bv \in \BBV\longmapsto
J(\bv)=(g(\bx),|\bv_{\sbtau}(\bx)|)_S$ is not differentiable at
zero. So, to bypass that hurdle we introduce the regularized
functional $J_\eps$ defined by
\begin{equation*}
v\in \BBV\mapsto J_\eps(v)=(g(\bx),
\sqrt{|v_{\sbtau}(\bx)|^2+\eps^2})_S~~,~~~ 0<\eps<<1~.
\end{equation*}
Clearly $J_\eps$ is convex and Gateaux differentiable with Gateaux
derivative $K_\eps$ defined on $\BBV$ and given by
\[
\langle K_\eps(\bu),\bv\rangle=\int_Sg\displaystyle
\frac{\bu_{\sbtau}\cdot\bv_{\sbtau}}{\sqrt{|\bu_{\sbtau}|^2+\eps^2}}\,ds\,.
\]
We briefly observe that $K_\eps$ is monotone. Indeed since $J_\eps$ is convex, for $\bu,\bv$ elements of $\BBV$ and $0<t<1$, $J_\epsilon(t\bu+(1-t)\bv)\leq tJ_\epsilon(\bu)+(1-t)J_\epsilon(\bv)$, which can be re-written as
\[
\displaystyle \frac{J_\epsilon(\bv+t(\bu-\bv))-J_\epsilon(\bv)}{t}\leq J_\epsilon(\bu)-J_\epsilon(\bv)\,.
\]
Then by taking the limit on both sides when $t$ goes to zero yields
\[
\langle K_\eps(\bv),\bu-\bv\rangle\leq J_\epsilon(\bu)-J_\epsilon(\bv)\,.
\]
Interchanging the role of $\bv$ and $\bu$, one gets instead
\[
\langle K_\eps(\bu),\bv-\bu\rangle\leq J_\epsilon(\bv)-J_\epsilon(\bu)\,.
\]
Putting together the later and former inequality, one has the
desired result
\begin{equation}
\langle K_\eps(\bu)-K_\eps(\bv),\bu-\bv\rangle\geq
0\,\,. \label{eq:monotone}
\end{equation}
The regularized form of (\ref{eq:2-4}) can be written as follows:
Find $\bu_\eps\in L^\infty(0,T;\BBV)$ satisfying (\ref{eq:2-3}) with
$\bu'_\eps\in L^2(0,T;\BBV')$ such that
for almost all $t$ and all $\bv\in \BBV$
\begin{eqnarray}
&&(\bu^{'}_{\eps}(t),\bv-\bu_\eps(t))+\nu
(\nabla\bu_\eps(t),\nabla(\bv-\bu_\eps(t))\,)+a\,(\bu_\eps(t),\bv-\bu_\eps(t))
\nonumber\\
&&+b\,(|\bu_\eps(t)|^\alpha\bu_\eps(t),\bv-\bu_\eps(t))
+J_\eps(\bv)-J_\eps(\bu_\eps(t))\geq (\fb(t),\bv-\bu_\eps(t))\,.\label{eq:2-5}
\end{eqnarray}
Since $J_\eps$ is differentiable, adopting the classical aruments in
\cite{DUVAUT}, one can state that (\ref{eq:2-5}) is equivalent to:
Find $\bu_\eps\in L^\infty(0,T;\BBV)$ satisfying (\ref{eq:2-3}) with
$\bu'_\eps\in L^2(0,T;\BBV')$ such that for almost all $t$
\begin{eqnarray}
&&(\bu^{'}_{\eps}(t),\bv)+\nu
(\nabla\bu_\eps(t),\nabla\bv)+a\,(\bu_\eps(t),\bv)+ b\,(|\bu_\eps(t)|^\alpha
\bu_\eps(t),\bv)\nonumber\\
&&+\langle K_\eps(\bu_\eps(t)),\bv\rangle= (\fb(t),\bv)~~\mbox{for all}~\bv\in \BBV\,\,.\label{eq:2-6}
\end{eqnarray}
Before proving the existence of a solution $\bu_\eps(t)$ of
(\ref{eq:2-6}), we first show how the pressure is constructed,
knowing the velocity. For that purpose, we begin by integrating
(\ref{eq:2-6}) on $[0,t]$, apply (\ref{eq:2-3}), and for $\bv\in
{\mathcal N}\cap \mathbf{L}^{\alpha+2}(\Omega)$; we introduce the
functional
\begin{eqnarray*}
{\mathcal H}(\bv)(t)&=&\int^t_0\Big[(\fb(s),\bv)-\nu
(\nabla\bu_\eps(s),\nabla\bv)-a\,(\bu_\eps(s),\bv)-
b\,(|\bu_\eps(s)|^\alpha \bu_\eps(s),\bv)\Big]~ds\\
&&-\int^t_0\langle
K_\eps(\bu_\eps(s)),\bv\rangle~ds-(\bu_\eps(t),\bv)+ (\bu_0,\bv)~,~~\text{for all}~~0\leq t\leq T~.
\end{eqnarray*}
One sees that ${\mathcal H}$ is linear and continuous on ${\mathcal N}\cap
\mathbf{L}^{\alpha+2}(\Omega)$, and according to (\ref{eq:2-6}) and
(\ref{eq:2-3}), it vanishes on $\BBV$. Now following
\cite{Raviart-Thomas}, [Theorem 2.5-1, page 54], for each $t\in
[0,T]$, there exists a unique function $\widetilde{p}_\eps(t)\in
L^2_0(\Omega)$ and a positive constant $C$ such that: for all
$\bv\in {\mathcal N}\cap \mathbf{L}^{\alpha+2}(\Omega)$,
\begin{eqnarray}
{\mathcal H}(\bv)(t)=(\Div\bv,\widetilde{p}_\eps(t))~,~~\label{eq:2-6b}\\
~~C\|\widetilde{p}_\eps(t)\|\leq
\displaystyle
\sup_{\sbv\in{\mathcal N}}\frac{(\Div\bv,\widetilde{p}_\eps(t))}{\|\bv\|_1}~~~.
\label{eq:2-7b}
\end{eqnarray}
Finally, we take the time derivative on both sides of
(\ref{eq:2-6b}); and we let
\begin{equation}
p_\eps(t)=\displaystyle \frac{d}{dt}\widetilde{p}_\eps(t)~,
\label{eq:2-8}
\end{equation}
in the resulting equation. Thus we have obtained the following
variational problem: Find
$\bu_\eps\in L^2(0,T;{\mathcal N}\cap \mathbf{L}^{\alpha+2}(\Omega))$,
$p_\eps\in L^2(0,T;L^2_0(\Omega))$ with
$\bu'_\eps\in L^2(0,T;({\mathcal N}\cap \mathbf{L}^{\alpha+2}(\Omega))')$
such that for almost all $t$ and all $q\in L^2(\Omega)~,~\bv\in {\mathcal N}\cap \mathbf{L}^{\alpha+2}(\Omega)$
\begin{eqnarray}
&&(\bu^{'}_{\eps}(t),\bv)+\nu
(\nabla\bu_\eps(t),\nabla\bv)+a\,(\bu_\eps(t),\bv)+ b\,(|\bu_\eps(t)|^\alpha
\bu_\eps(t),\bv)\nonumber\\
&&-(\Div\bv,p_\eps(t))+\langle K_\eps(\bu_\eps(t)),\bv\rangle=(\fb(t),\bv)\,,
\label{eq:2-9}\\
&&(\Div\bu_\eps(t),q)=0\,\,\nonumber\\
&&\bu_\eps(0)=\bu_0~.
\nonumber
\end{eqnarray}
It is clear that the variational problems (\ref{eq:2-6}) and
(\ref{eq:2-9}) are equivalent, with the regularized pressure
described by (\ref{eq:2-6b}), (\ref{eq:2-7b}) and
(\ref{eq:2-8})\,.\newline

Step 2: {\it Faedo-Galerkin approximation}.\newline
We let
\begin{eqnarray*}
\BBH&=&\{\bv\in \mathbf{L}^2(\Omega)~,~~\Div\bv=0~,~~v_n|_{\partial
\Omega}=0\}\cap \mathbf{L}^{\alpha+2}(\Omega)~.
\end{eqnarray*}
One readily observes that $\BBV$ is compactly embedded in $\BBH$\,. For the slip
boundary condition, we introduce the Stokes operator defined on a
subspace of $\BBV$ constructed in \cite{Heywood} as follows; for every $\fb\in \BBH$, there exists a unique
$\bv\in\BBV$ such that
\begin{equation}
\label{ONE}
(\nabla \bv, \nabla \phi)=(\fb,\phi), \forall \phi \in \BBV\,.
\end{equation}
Moreover, for every $\bv\in \BBV$, there exists a unique $\fb\in
\BBH$ such that \eqref{ONE} holds. Then \eqref{ONE} defines a
one-to-one mapping between $\fb\in \BBH$ and $\bv\in D(A)$, where
$D(A)$ is a subspace of $\BBV$. Hence, $A\bv=\fb$ defines the Stokes
operator $A:D(A)\rightarrow \BBH$. Its inverse $A^{-1}$ is compact
and self-adjoint as a mapping from $\BBH$ to $\BBH$ and possesses an
orthogonal sequence of eigenfunctions $\bpsi_k$ which are complete
in $\BBH$ and $\BBV$;
\begin{equation}\label{TWO}
A\bpsi_k=\lambda_k \bpsi_k.
\end{equation}
Let $\BBV_m$ be the subspace of $\BBV$ spanned by
$\psi_1,\cdots,\psi_m$, that is
\[
\BBV_m=\{\psi_1,\psi_2,\psi_3,\cdots,\psi_m\}~.
\]
We consider the following ordinary differential equation:
Find~$\bu_{\eps,m}(t)\in \BBV_m$~such that for all~$\bv\in \BBV_m$;
\begin{equation}
\label{eq:2-10}
\begin{split}
(\bu^{'}_{\eps,m}(t)~,~\bv)+\nu
(\nabla\bu_{\eps,m}(t)~,~\nabla\bv)+
a\,(\bu_{\eps,m}(t)~,~\bv)\\+b\,(|\bu_{\eps,m}(t)|^\alpha
\bu_{\eps,m}(t)~,~\bv)
+\langle K_\eps(\bu_{\eps,m}(t))~,~\bv\rangle=
(\fb(t),\bv)~~,\\
\bu_{\eps,m}(0)\rightarrow \bu_\eps(0)=\bu_0 \in \BBV_m~.
\end{split}
\end{equation}
As far as the existence of $\bu_{\eps,m}(t)$ defined by (\ref{eq:2-10}) is concerned, we note that the mapping
\[
{\mathcal K}:\bw\longmapsto (\fb,\bv)- \nu (\nabla\bw~,~\nabla\bv)-
a\,(\bw~,~\bv)-b\,(|\bw|^\alpha \bw~,~\bv)-\langle
K_\eps(\bw)~,~\bv\rangle~,
\]
is locally Lipschitz thanks to the nature of the operators involved.
It then follows from the theory of ordinary differential equations
that (\ref{eq:2-10}) has a solution $\bu_{\eps,m}$ defined on
$[0,t_{\eps,m}]$, $t_{\eps,m}>0$. Hereafter, $C$ denotes a constant
independent of $m$, and depending only on the data such as $\Omega$,
and whose value may be different in each inequality. Next, we derive
some {\it a priori} estimates and deduce that $t_{\eps,m}$ does not
depend on $\eps$ or $m$. Concerning the later property, is should be
mentioned from \cite{DUVAUT,LIONS}, that it suffice to derive {\it a
priori} estimates of the solution with the right hand side
independent of $m$ and $\eps$.

{\it Step 3: Some a priori estimates}.\newline
First we let $\bv=\bu_{\eps,m}(t)$ in (\ref{eq:2-10}). After using Young's inequality, one obtains
\begin{eqnarray}
&&\displaystyle \frac{d}{dt}\|\bu_{\eps,m}(t)\|^2+2\nu \|\nabla
\bu_{\eps,m}(t)\|^2+a\|\bu_{\eps,m}(t)\|^2
+2b\|\bu_{\eps,m}(t)\|^{\alpha+2}_{L^{\alpha+2}}\nonumber\\
&&+2
\langle K_\eps(\bu_{\eps,m}(t))~,~\bu_{\eps,m}(t)\rangle
\leq
\displaystyle\frac{||f(t)||^2}{a}~,
\label{eq:2-11a}
\end{eqnarray}
which by integrating over $[0,T^\sharp]$ for $T^\sharp\leq t_{\eps,m}$, and using (\ref{eq:monotone}), yields
\begin{eqnarray}
&&\sup_{0\le t\le T^\sharp}||\bu_{\eps,m}(t)||^2+2\nu \int_{0}^{T^\sharp}||\nabla
\bu_{\eps,m}(t)||^2 dt+a\int_{0}^{T^\sharp}||\bu_{\eps,m}(t)||^2dt\nonumber\\
&&+2b\int^{T^\sharp}_{0}\|\bu_{\eps,m}(t)\|^{\alpha+2}_{L^{\alpha+2}} dt\leq \displaystyle \frac{1}{a}\int_{0}^{T^\sharp}||\fb(t)||^2 dt+||\bu_0||^2<\infty\,
\label{eq:2-11}
\end{eqnarray}
since by assumption $\fb\in L^2(Q)$\,.
Now let $\bv=\bu^{'}_{\eps,m}(t)$ in (\ref{eq:2-10}). For $0\leq t\leq T^\sharp$, Young's inequality yields
\begin{eqnarray*}
&&||\bu^{'}_{\eps,m}(t)||^2+\frac{d}{dt}\left[ \nu
||\nabla \bu_{\eps,m}(t)||^2 +
a||\bu_{\eps,m}(t)||^2+\frac{2b}{\alpha+2}\|
\bu_{\eps,m}(t)\|^{\alpha+2}_{L^{\alpha+2}}\right]\\ && +\frac{d}{dt}\left[2K_\eps(\bu_{\eps,m}(t))\right]\leq ||\fb(t)||^2~,
\end{eqnarray*}
which leads to
\begin{eqnarray}
&&\int_{0}^{T^\sharp}||\bu^{'}_{\eps,m}(t)||^2dt+\nu||\nabla \bu_{\eps,m}(t)||^2 + a||\bu_{\eps,m}(t)||^2+\frac{2b}{\alpha+2}\|\bu_{\eps,m}(t)\|^{\alpha+2}_{L^{\alpha+2}} \nonumber\\ 
&& +2K_\eps(\bu_{\eps,m}(t))\leq ||\fb||^{2}_{\dis L^2\dis(0,T^\sharp;L^2)}+\Phi(0),\label{eq:2-12}
\end{eqnarray}
where
\begin{equation}
\Phi(0)=\nu||\nabla \bu_0||^2 +
a||\bu_0||^2+\frac{2b}{\alpha+2}\|\bu_0\|^{\alpha+2}_{L^{\alpha+2}}+2\int_S
g\sqrt{|\bu_0|^2+1}~ds~. \nonumber
\end{equation}
It is manifest that the right hand sides of the {\it a priori}
estimates obtained in (\ref{eq:2-11}) and (\ref{eq:2-12}) are
independent of $m$ and $\eps$. We then conclude that $t_{\eps,m}$ is
independent of $\eps$ and $m$ following the arguments discussed in
length by \cite{DUVAUT,LIONS}.

Step 4: {\it Passage to the limit}.\newline
We need to pass to the limit when $m$ approaches infinity and $\eps$ approaches zero. We start by fixing
$\eps$ and study the sequence
$m\longmapsto\bu_{\eps,m}$.\newline Based on (\ref{eq:2-11}) and
(\ref{eq:2-12}), it is clear that when $m\rightarrow \infty$,
\begin{align}
&\bu_{\eps,m}~\text{remains bounded in}~L^\infty(0,T;\BBH),~\nonumber\\
&|\bu_{\eps,m}|^\alpha\bu_{\eps,m}~\text{remains bounded in}~L^{\frac{\alpha+2}{\alpha+1}}(0,T;\mathbf{L}^{\frac{\alpha+2}{\alpha+1}}(\Omega))~~,
\label{eq:boundedset}\\
&\bu^{'}_{\eps,m}~\text{remains bounded in}~L^2(0,T;\mathbf{L}^2(\Omega))~. \nonumber
\end{align}
From a consequence of the result of {\it Dunford-Pettis} \cite{Yosida},
it is possible to extract from $(\bu_{\eps,m})_m$ a
subsequence, denoted again by $(\bu_{\eps,m})_m$ such that
\begin{align}
&\bu_{\eps,m}\longrightarrow \bu_\eps \text{ weak star in }~ L^\infty(0,T;\BBH)\label{eq:2-13}\\
&\bu_{\eps,m}\longrightarrow \bu_\eps \text{ weak star in }~ L^\infty(0,T;\BBV_m)\label{eq:2-14}\\
&|\bu_{\eps,m}|^\alpha\bu_{\eps,m}\longrightarrow \chi_\eps~ \text{weak star in}~
L^{\frac{\alpha+2}{\alpha+1}}
\left(0,T;\mathbf{L}^{\frac{\alpha+2}{\alpha+1}}(\Omega)\right)\label{eq:2-15}\\
&\bu^{'}_{\eps,m}\longrightarrow \bu^{'}_\eps \text{
weak in }~ L^2(0,T;\BBH).\label{eq:2-16}
\end{align}
The convergence results (\ref{eq:2-13}), and (\ref{eq:2-16}) imply in particular that \begin{align}
&\bu_{\eps,m}~\text{remains in a bounded set of}~H^1(Q)\,.
\label{eq:2-16b}
\end{align}
But from Rellich-Kondrachoff, the embedding $H^1(Q)\longmapsto L^2(Q)$ is compact. So one can extract from $(\bu_{\eps,m})$ a
subsequence, denoted again by $(\bu_{\eps,m})$ such that
\begin{align}
&\bu_{\eps,m}\longrightarrow \bu_\eps~\text{strong in}~L^2(0,T;\BBH)~\text{and a.e.~in}~Q\,.
\label{eq:2-17}
\end{align}
Next, it follows from (\ref{eq:2-15}) and (\ref{eq:2-17}) and Lemma
1.3 in \cite{LIONS} (page 12) that $\chi_\eps=|\bu_\eps|^\alpha\bu_\eps$~.\newline
It remains to be shown that
\begin{equation}
\label{eq:2-19}
K_\eps(\bu_{\eps,m})\longrightarrow K_\eps(\bu_\eps) \text{ weak star in } L^\infty(0,T,\BBV_m')~.
\end{equation}
Firstly from (\ref{eq:2-14})
\begin{equation}\label{eq:2-20}
K_\eps(\bu_{\eps,m})\longrightarrow \beta_\eps \text{ weak star in } L^\infty(0,T,\BBV_m')~.
\end{equation}
Passing to the limit in (\ref{eq:2-10}), one obtains
\begin{eqnarray}
&&(\bu'_\eps(t),\bv)+\nu
(\nabla\bu_\eps(t)~,~\nabla\bv)+a\,(\bu_\eps(t)~,~\bv)+b\,(|\bu_\eps(t)|^\alpha
\bu_\eps(t)~,~\bv)\nonumber\\
&&+\langle\beta_\eps,\bv\rangle=(\fb(t),\bv), \forall
\bv\in \BBV_m~.\label{eq:2-21}
\end{eqnarray}
For any $\bw \in L^1(0,T;\BBV_m)$, since $K_\eps(\cdot)$ is
monotone (see \ref{eq:monotone}),
\[
\langle K_\eps(\bu_{\eps,m}(t))~,~\bu_{\eps,m}(t)\rangle\geq \langle
K_\eps(\bu_{\eps,m}(t))~,~\bw\rangle+\langle
K_\eps(\bw)~,~\bu_{\eps,m}(t)-\bw\rangle~~,
\]
but from (\ref{eq:2-10})
\begin{eqnarray*}
\langle
K_\eps(\bu_{\eps,m}(t))~,~\bu_{\eps,m}(t)\rangle&=&(\fb(t),\bu_{\eps,m}(t))-
(\bu^{'}_{\eps,m}(t)~,~\bu_{\eps,m}(t))\\&& - \nu
(\nabla\bu_{\eps,m}(t)~,~\nabla\bu_{\eps,m}(t))-a\,(\bu_{\eps,m}(t)~,~\bu_{\eps,m}(t))\\
&& -b\,(|\bu_{\eps,m}(t)|^\alpha \bu_{\eps,m}(t)~,~\bu_{\eps,m}(t))~.
\end{eqnarray*}
The former and latter equations give
\begin{eqnarray*}
(\fb(t),\bu_{\eps,m}(t))-\frac{1}{2}\frac{d}{dt}\|\bu_{\eps,m}(t)\|^2-\nu
\|\nabla\bu_{\eps,m}(t)\|^2-a\|\bu_{\eps,m}(t)\|^2
-b\|\bu_{\eps,m}(t)\|^{\alpha+2}_{L^{\alpha+2}}\\
\geq \langle K_\eps(\bu_{\eps,m}(t))~,~\bw\rangle+\langle
K_\eps(\bw)~,~\bu_{\eps,m}(t)-\bw\rangle~.
\end{eqnarray*}
So, integrating with respect to $t$ on $[0,T]$, yields
\begin{eqnarray}
&&\int^T_0(\fb(t),\bu_{\eps,m}(t))dt- \displaystyle
\frac{1}{2}\|\bu_{\eps,m}(T)\|^2+\displaystyle
\frac{1}{2}\|\bu_{\eps,m}(0)\|^2\nonumber\\
&&-\int^T_0\displaystyle \left[\nu
\|\nabla\bu_{\eps,m}(t)\|^2+a\|\bu_{\eps,m}(t)\|^2+b\|\bu_{\eps,m}(t)\|^{\alpha+2}_{L^{\alpha+2}}\right]dt\label{eq:2-22}\\
&&\geq \int^T_0\left[\langle
K_\eps(\bu_{\eps,m}(t))~,~\bw(t)\rangle+\langle
K_\eps(\bw(t))~,~\bu_{\eps,m}(t)-\bw(t)\rangle\right]dt~.\nonumber
\end{eqnarray}
Next, we take $\bv=\bu_{\eps,m}(t)$ in (\ref{eq:2-21}), and combined
the resulting equation with (\ref{eq:2-22}), which yields (after
taking the limit as $m$ approaches to infinity)
\begin{eqnarray}
\int^T_0\langle
\beta_\eps-K_\eps(\bw(t))~,~\bu_\eps(t)-\bw(t)\rangle
dt\geq0~.\label{eq:2-23}
\end{eqnarray}
At this juncture, we let $\bu_\eps(t)-\bw(t)=\pm\bq$ with $\bq\in
L^2(0,T;\BBV_m)$\,. Thus (\ref{eq:2-23}) leads to
\begin{eqnarray*}
\int^T_0\langle \beta_\eps-K_\eps(\bw(t))~,~\bq\rangle dt=0~~,~
\end{eqnarray*}
from which we deduce the desired convergence result (\ref{eq:2-19}).
We have established that as $m$ goes to infinity, the sequence
$(\bu_{\eps,m}(t))_m$ converges to $\bu_\eps(t)$ in some sense with
$\bu_\eps(t)$, the solution of
\begin{equation}
\begin{split}
\label{eq:2-24}
(\bu'_\eps(t),\bv)+\nu(\nabla\bu_\eps(t)~,~\nabla\bv)+a\,(\bu_\eps(t)~,~\bv)
+b\,(|\bu_\eps(t)|^\alpha \bu_\eps(t)~,~\bv)\\
+\langle K_\eps(\bu_\eps(t)),\bv\rangle=(\fb(t),\bv)~~,~~\text{for all}~~~
\bv\in \BBV_m~.
\end{split}
\end{equation}

Since $\cup_m \BBV_m$ is dense in $\BBV$, we can conclude that
(\ref{eq:2-24}) holds true for $\bv$ in $\BBV$. Therefore, we have
established that there exists a function $\bu_\eps$ uniformly
bounded with respect to $\eps$ in $L^\infty(0,T, \BBH\cap\BBV\cap
\mathbf{L}^{\alpha+2}(\Omega))$ such that $\bu'_\eps$ is uniformly
bounded with respect to $\eps$ in $L^2(0,T, \BBH)$ and $\bu_\eps$
satisfies (\ref{eq:2-24}).\newline
Our final task in the paragraph is to consider the limit as
$\eps$ goes to zero.\newline
First, we take the limit on both sides of (\ref{eq:2-11}) and
(\ref{eq:2-12}), one has
\begin{eqnarray}
&&\sup_{0\le t\le T}||\bu_\eps(t)||^2+2\nu \int_0^T||\nabla
\bu_\eps(t)||^2 dt+a\int_0^T||\bu_\eps(t)||^2dt\nonumber\\
&&+2b\int^T_0\|\bu_\eps(t)\|^{\alpha+2}_{L^{\alpha+2}} dt\le
\displaystyle \frac{1}{a}\int_0^T||\fb(t)||^2
dt+||\bu_0||^2,~~\label{eq:2-25}
\end{eqnarray}
and 
\begin{eqnarray}
&&\int_0^T||\bu'_\eps(t)||^2dt+\nu||\nabla \bu_\eps(t)||^2
+a||\bu_\eps(t)||^2+\frac{2b}{\alpha+2}\|\bu_\eps(t)\|^{\alpha+2}_{L^{\alpha+2}}
\nonumber\\
&&\le \int_0^T||\fb(t)||^2dt+\Phi(0)~.\label{eq:2-26}
\end{eqnarray}
Thus we can extract from
$(\bu_\eps)_\eps$ a subsequence still denoted by $(\bu_\eps)_\eps$
such that
\begin{align}
\bu_\eps\longrightarrow \bu \text{ weak star in } L^\infty(0,T,\BBH)\label{eq:2-27}\\
\bu_\eps\longrightarrow \bu \text{ weak star in } L^\infty(0,T,\BBV)\label{eq:2-28}\\
|\bu_\eps|^\alpha\bu_\eps\longrightarrow \chi \text{ weak star in } L^{\frac{\alpha+2}{\alpha+1}}
\left(0,T,\mathbf{L}^{\frac{\alpha+2}{\alpha+1}}(\Omega)\right)\label{eq:2-29}\\
\bu'_\eps \longrightarrow \bu' \text{ weak in } L^2(0,T,\BBH)~.
\end{align}
Arguing as before we can prove that
\begin{align}
\bu_\eps\longrightarrow \bu \text{ strong in } L(0,T;\BBH) \text{ and a.e. in}\,\,Q,\label{eq:2-30}\\
|\bu_\eps|^\alpha \bu_\eps\longrightarrow |\bu|^\alpha \bu \text{
weak in } L^{\frac{\alpha+2}{\alpha+1}}\left(0,T;\mathbf{L}^{\frac{\alpha+2}{\alpha+1}}(\Omega)\right)~.
\end{align}
\noindent  Let $\bv\in L^2(0,T,\BBV)$, from (\ref{eq:2-24}), it
follows that
\begin{eqnarray}
&&(\bu'_\eps(t),\bv-\bu_\eps(t))+\nu
(\nabla\bu_\eps(t)~,~\nabla(\bv-\bu_\eps(t)))+a\,(\bu_\eps(t)~,~\bv-\bu_\eps(t))\nonumber\\
&&+b\,(|\bu_\eps(t)|^\alpha
\bu_\eps(t)~,~\bv-\bu_\eps(t))+J_\eps(\bv)-J_\eps(\bu_\eps(t))\nonumber\\
&&=(\fb(t),\bv-\bu_\eps(t))+J_\eps(\bv)-J_\eps(\bu_\eps(t))-\langle
K_\eps(\bu_\eps(t)),\bv-\bu_\eps(t)\rangle~.\label{eq:2-31}
\end{eqnarray}
Integrating (\ref{eq:2-31}) with respect to $t$ along $[0,T]$ and
taking into account the fact that
$J_\eps(\bv)-J_\eps(\bu_\eps(t))-\langle
K_\eps(\bu_\eps(t)),\bv-\bu_\eps(t)\rangle\geq 0$, one obtains
\begin{equation}
\label{eq:2-32}
\begin{split}
\int_0^T\left(
(\bu'_\eps(t),\bv)+\nu(\nabla\bu_\eps(t),\nabla\bv)+
a(\bu_\eps(t),\bv)+b(|\bu_\eps(t)|^\alpha\bu_\eps(t),\bv)\right)dt\\ +\int_0^T \left (J_\eps(\bv)-
(\fb(t),\bv-\bu_\eps(t))\right)dt\\
\ge
\frac{1}{2}||\bu_\eps(T)||^2-\frac{1}{2}||\bu_{\eps0}||^2+\int_0^T\left(
a||\bu_\eps(t)||^2+b\int_\Omega
|\bu_\eps(t)|^{\alpha+2}dx\right)dt\\+\int_0^tJ_\eps(\bu_\eps(t))dt~.
\end{split}
\end{equation}
Since $\bu_\eps\longrightarrow \bu$ weak star in
$L^\infty(0,T,\BBV)$, and $J_\eps$ is a convex and continuous functional on $\BBV$, one has
\begin{equation}
\label{eq:2-33} \displaystyle \lim\inf_{\eps
\rightarrow0}\int_0^TJ_\eps(\bu_\eps(t))dt\ge \int_0^TJ(\bu(t))dt~.
\end{equation}
By using (\ref{eq:2-33}), we infer
from (\ref{eq:2-32}) that
\begin{equation*}
\begin{split}
\int_0^T\big((
\bu'(t),\bv)+\nu(\nabla\bu(t),\nabla\bv)+a(\bu(t),\bv)+b(|\bu(t)|^\alpha
\bu(t),\bv)\\+
J(\bv)-(\fb(t),\bv-\bu(t))\big)dt\\
\ge \frac{1}{2}||\bu(T)||^2-\frac{1}{2}||\bu_0||^2+\int_0^T\left(
a||\bu(t)||^2+b\int_\Omega
|\bu(t)|^{\alpha+2}dx\right)dt+\int_0^TJ(\bu(t))dt\\
= \int^T_0\left[(\bu'(t),\bu(t))+
a\,(\bu(t),\bu(t))+b\,(|\bu(t)|^\alpha\bu(t),\bu(t))+J(\bu(t))\right]dt\\
\end{split}
\end{equation*}
which by arguing as in \cite{DUVAUT}, pages 56-57, yields
\begin{equation*}
\begin{split}
(\bu'(t),\bv-\bu(t))+\nu
(\nabla\bu(t),\nabla(\bv-\bu(t)))+a(\bu(t),\bv-\bu(t))\\+b(|\bu(t)|^\alpha
\bu(t),\bv-\bu(t))
+J(\bv)-J(\bu(t))\geq (\fb(t),\bv-\bu(t))~~~\text{for all}~~\bv\in \BBV~.
\end{split}
\end{equation*}
We then conclude that
\begin{theorem}
The variational problem (\ref{eq:2-6}) admits at least a weak
solution, which moreover satisfies;
\begin{eqnarray}
&&\displaystyle \sup_{0\leq t\leq T}\|\nabla \bu(t)\|\leq C~~,~~
\displaystyle \int^T_0\|\bu'(t)\|^2dt\leq C~,
\end{eqnarray}
where $C$ is a positive constant depending on the data.
\label{existence}
\end{theorem}
Having obtained the velocity, we shall indicate how the pressure is
constructed. First, we recall that from $(\ref{eq:2-9})_1$,
\begin{equation*}
\begin{split}
(\Div\bv,p_\eps(t))=(\bu'_\eps(t),\bv)+\nu
(\nabla\bu_\eps(t),\nabla\bv)+a\,(\bu_\eps(t),\bv)\\+ b\,(|\bu_\eps(t)|^\alpha
\bu_\eps(t),\bv)
+\langle K_\eps(\bu_\eps(t)),\bv\rangle-(\fb(t),\bv)\,,
\end{split}
\end{equation*}
but since $p_\eps(t)\in L^2_0(\Omega)$, following
\cite{Girault-Raviart}, one can find a positive constant $C$ such
that
\begin{eqnarray*}
C\|p_\eps(t)\|\leq \displaystyle \sup_{\sbv\in {\mathcal
N}}\frac{(\Div\bv,p_\eps(t))}{\|\bv\|_1}~~.
\end{eqnarray*}
Now, combining the former and latter equations and the continuity of operators involved, one obtains
\begin{equation*}
\begin{split}
C\|p_\eps(t)\|\leq 
\|\bu'_\eps(t)\|+\nu \|\nabla\bu_\eps(t)\|+a\|\bu_\eps(t)\|+ b\|\bu_\eps(t)\|^{\alpha+1}_{L^{2\alpha+2}}\\+\|K_\eps(\bu_\eps(t))\|_{{\mathcal V}'}+\|\fb(t)\|
\end{split}
\end{equation*}
Equivalently,
\begin{equation*}
\begin{split}
C\|p_\eps(t)\|&\leq \|\bu'_\eps(t)\|+\nu \|\nabla\bu_\eps(t)\|+a\|\bu_\eps(t)\|+
C(b,\Omega,\alpha)\|\bu_\eps(t)\|^{\alpha+1}_{L^6}\\&\quad +C(\Omega)\|g\|_{L^\infty(S)}\|\bu_\eps(t)\|_1+\|\fb(t)\|,\\
&\leq \|\bu'_\eps(t)\|+\nu \|\nabla\bu_\eps(t)\|+a\|\bu_\eps(t)\|+
C(b,\Omega,\alpha)\|\nabla\bu_\eps(t)\|^{\alpha+1}\\ &
+C(\Omega)\|g\|_{L^\infty(S)}\|\bu_\eps(t)\|_1+\|\fb(t)\|~
\end{split}
\end{equation*}
which by Young's inequality and integrating the resulting inequality over
$[0,T]$, yields (after utilization of (\ref{eq:2-25}) and
(\ref{eq:2-26}))
\begin{equation}
\begin{split}
\int^T_0\|p_\eps(t)\|^2dt\leq
C\int^T_0\|\bu'_\eps(t)\|^2+C\int^T_0\|\nabla\bu_\eps(t)\|^2dt+C\int^T_0\|\bu_\eps(t)\|^2\\ +
C\int^T_0\|\nabla\bu_\eps(t)\|^{2\alpha+2}dt
+C\|g\|^2_{L^\infty(S)}\int^T_0\|\bu_\eps(t)\|^2_1\\+C\int^T_0\|\fb(t)\|^2dt<\infty~~,
\label{eq:2-34}
\end{split}
\end{equation}
$C$ being a positive constant depending on the parameters and the
domain of the problem. Then we can select from $p_\eps(t)$ a sequence,
again denoted by $p_\eps(t)$, such that
\begin{eqnarray}
p_\eps\longrightarrow p~~\text{weakly
in}~~~L^2(0,T;L^2_0(\Omega))~~. \label{eq:2-35}
\end{eqnarray}
Next, one observes that (\ref{eq:2-9}) can be re-written as
\begin{equation*}
\begin{split}
(\bu'_\eps(t),\bv-\bu_\eps(t))+\nu
(\nabla\bu_\eps(t),\nabla(\bv-\bu_\eps(t)))+a\,(\bu_\eps(t),\bv-\bu_\eps(t))\\
+
b\,(|\bu_\eps(t)|^\alpha
\bu_\eps(t),\bv-\bu_\eps(t))
-(\Div(\bv-\bu_\eps(t))~,~p_\eps(t))\\+J_\eps(\bv)-J_\eps(\bu_\eps(t))-
(\fb(t),\bv-\bu_\eps(t))\geq
0~~\,\text{for all}~~\bv\in {\mathcal N}\cap \mathbf{L}^{\alpha+2}(\Omega)\,,\\
(\Div\bu_\eps(t),q)=0\,\,,
\text{for all}~~q\in L^2(\Omega)\,,
\end{split}
\end{equation*}
which by integration over the time interval $[0,T]$ and passage to
the limit (as $\eps\rightarrow 0$ ) yields, (after utilization of
the identity $(\Div\bu_\eps(t),q)=0$ for all $q\in L^2(\Omega)$)
\begin{equation*}
\begin{split}
\int^T_0\big[(\bu'(t),\bv-\bu(t))+\nu
(\nabla\bu(t),\nabla(\bv-\bu(t)))+a\,(\bu(t),\bv-\bu(t))\big]dt\\
+\int^T_0\big[-(\Div(\bv-\bu(t))~,~p(t))+J(\bv)-J(\bu(t))-(\fb(t),\bv-\bu(t))\big]dt
\\ + \int_0^T b\,(|\bu(t)|^\alpha
\bu(t),\bv-\bu(t)) dt\geq
0,
\end{split}
\end{equation*}
$\text{for all}~~\bv\in {\mathcal N}\cap \mathbf{L}^{\alpha+2}(\Omega)$. Also, 
$(\Div\bu(t),q)=0\,\,
\text{for all}~~q\in L^2(\Omega)~.$

Finally, arguing as in \cite{DUVAUT}, pages 56-57, one obtains
\begin{equation}
\begin{split}
(\bu'(t),\bv-\bu(t))+\nu
(\nabla\bu(t),\nabla(\bv-\bu(t)))+a\,(\bu(t),\bv-\bu(t))\\+ b\,(|\bu(t)|^\alpha
\bu(t),\bv-\bu(t))
-(\Div(\bv-\bu(t))~,~p(t))+J(\bv)-J(\bu(t))\\
\geq (\fb(t),\bv-\bu(t))
\end{split}
\end{equation}
\text{for all} $
\bv\in {\mathcal N}\cap \mathbf{L}^{\alpha+2}(\Omega)$. Moreover,
$(\Div\bu(t),q)=0\,\,\,\text{ for all }~q\in L^2(\Omega).$
\section{Continuous dependence on the data}

In this section, our focus is to establish some qualitative
properties of the weak  solutions  in Theorem \ref{existence}. In
particular, we show that the solutions depend continuously on
initial velocity, external force as well as the Forchheimer's and Brinkman's
coefficients. We recall that such results in the literature are
sometimes referred to as structural stability.
\newline
We first claim that

\begin{theorem}\label{CONTINUOUS1}
Let $\bu_i$ be the solution of (\ref{eq:2-3}) with respect to
$\bu_{i0}, \fb_i$, $i=1,2$. Then there exists a positive constant
$C$, depending on $a,\nu$ and $\Omega$ such that
\begin{equation}
||\bu_1(t)-\bu_2(t)||^2\leq
e^{-C\,t}||\bu_1(0)-\bu_2(0)||^2+\displaystyle
\int^t_0e^{C(-t+s)}||\fb_2(s)-\fb_1(s)||^2ds\,.
\end{equation}
\end{theorem}

This theorem implies in particular the following uniqueness result.
\begin{corollary}
The problem (\ref{eq:2-3}) has one and only one solution.
\end{corollary}
{\bf Proof theorem \ref{CONTINUOUS1}}.\quad The functions $\bu_1$
and $\bu_2$ satisfy respectively:
\begin{equation}
\label{eq:3-1}
\begin{split}
(\partial_t \bu_1,\bv-\bu_1)+\nu
(\nabla\bu_1,\nabla(\bv-\bu_1))+a\,(\bu_1,\bv-\bu_1)+b\,(|\bu_1|^\alpha
\bu_1,\bv-\bu_1)\\
+J(\bv)-J(\bu_1)\ge (\fb_1,\bv-\bu_1)~~\mbox{for
all}~~\bv\in \BBV~.
\end{split}
\end{equation}
and
\begin{equation}\label{eq:3-2}
\begin{split}
(\partial_t\bu_2,\bv-\bu_2)+\nu
(\nabla\bu_2,\nabla(\bv-\bu_2))+a\,(\bu_2,\bv-\bu_2)+b\,(|\bu_2|^\alpha
\bu_2,\bv-\bu_2)\\+J(\bv)-J(\bu_2)\ge (\fb_2,\bv-\bu_2)~~\mbox{for
all}~~\bv\in \BBV~.
 \end{split}
\end{equation}
Setting $\bv=\bu_2$ in (\ref{eq:3-1}) and $\bv=\bu_1$ in
(\ref{eq:3-2}) and adding the resulting inequalities, it follows that
\begin{equation*}
\begin{split}
\displaystyle
\frac{1}{2}\frac{d}{dt} ||\bw(t)||^2+\nu||\nabla \bw||^2 +a ||\bw(t)||^2\\
+b(|\bu_2|^\alpha \bu_2-|\bu_1|^\alpha \bu_1,\bw(t)) \le
(\fb_2-\fb_1,\bw(t))~,
\end{split}
\end{equation*}
where $\bw(t)=\bu_2(t)-\bu_1(t)$ and $\bw_0=\bu_{20}-\bu_{10}$.
\noindent Since $T(\zeta)=|\zeta|^\alpha \zeta$ is monotone then
\begin{equation*}
(|\bu_2|^\alpha \bu_2-|\bu_1|^\alpha \bu_1,\bw(t))\ge 0~.
\end{equation*}
Therefore
\begin{equation}\label{eq:3-3}
\displaystyle \frac{d}{dt} ||\bw(t)||^2+C(\nu,a,\Omega)
||\bw(t)||^2\leq C(a,\Omega)||\fb_2-\fb_1||^2~~,
\end{equation}
where Poincar\'{e}'s inequality has been used. We readily deduce the
desired result from (\ref{eq:3-3}) using Gronwall's lemma\,.
\QED\newline

In line of theorem \ref{CONTINUOUS1}, one can state the following
result.
\begin{theorem}
The weak solutions of problem (\ref{eq:2-3}) constructed in theorem
\ref{existence} depends continuously with respect to the $L^2$ norm
on:
\begin{enumerate}
\item[(a)]\quad the Forchheimer coefficient $b$, and
\item[(b)]\quad the Brinkman coefficient $\nu$\,.
\end{enumerate}
\end{theorem}
The proof follows mutatis mutandis the proof of theorem
\ref{CONTINUOUS1}\,.

\section{Stability of stationary solutions}

Hereafter, we study the stability of stationary solutions to
(\ref{eq:2-3}).\newline We assume that the apply force $\fb$ is
independent of time, and we consider the following stationary
problem
\begin{equation}
\label{eq:4-1}
\begin{cases}
-\nu \Delta \bu+a\bu +b|\bu|^\alpha \bu-\nabla p=\fb, \text{ in } \Omega\,,\\
\Div \bu=0, \text{ in } \Omega\,,\\
\bu=0 \text{ on } \Gamma,\\
\bu\cdot\bn=0~~,~~\mbox{and}~~-\bsigma_{\sbtau}\in
g\partial|\bu_{\sbtau}| \text{ on } S\,.
\end{cases}
\end{equation}
Here, we always assume that $\alpha \in [1,2]$, $\gamma, a,b>0$.

It is clear that the velocity satisfies the simpler variational
problem
\begin{equation}
\label{eq:4-2}
\begin{cases}
\mbox{Find}~~\bu\in \BBV~\mbox{such that for all}~\bv\in \BBV~,\\
\nu(\nabla\bu ,\nabla(\bv-\bu))+a\,(\bu, \bv-\bu)+b(|\bu|^\alpha\bu,
\bv-\bu)+J(\bv)-J(\bu)\ge (\fb,\bv-\bu)~.
 \end{cases}
\end{equation}
It can be shown as in \cite{DUVAUT} that there exists a unique
$\widetilde{\bu}\in \BBV$ such that (\ref{eq:4-2}) holds true, and
one has the following

\begin{theorem}
The weak solution $\bu$ of (\ref{eq:2-3}) constructed in theorem
\ref{existence} converges to the unique solution $\widetilde{\bu}$
to (\ref{eq:4-2}) exponentially as $t$ goes to infinity. More
precisely, we have the following estimate
\begin{equation}
 ||\bu(t)-\widetilde{\bu}||^2\le ||\bu_0-\widetilde{\bu}||^2e^{-2(a+\nu)t}~~,~~\mbox{for all}~t\geq 0~.
\end{equation}
\end{theorem}
{\bf proof.}\quad We let $\bv=\bu(t)$ in (\ref{eq:4-2}), thus
\begin{equation*}
\begin{split}
\nu(\nabla\widetilde{\bu}
,\nabla(\bu(t)-\widetilde{\bu}))+a\,(\widetilde{\bu},
\bu(t)-\widetilde{\bu})+b(|\widetilde{\bu}|^\alpha\widetilde{\bu},
\bu(t)-\widetilde{\bu})\\+J(\bu(t))-J(\widetilde{\bu})\ge
(\fb,\bu(t)-\widetilde{\bu})~.
\end{split}
\end{equation*}
Next, for $\bv=\widetilde{\bu}$ in (\ref{eq:2-3}), one has
\begin{equation*}
\begin{split}
(\bu'(t),\widetilde{\bu}-\bu(t))+\nu
(\nabla\bu(t),\nabla(\widetilde{\bu}-\bu(t))\,)+a\,(\bu(t),\widetilde{\bu}-\bu(t))\\+
b\,(|\bu(t)|^\alpha
\bu(t),\widetilde{\bu}-\bu(t))+J(\widetilde{\bu})-J(\bu(t))\ge
(\fb,\widetilde{\bu}-\bu(t))\,\,\,.
\end{split}
\end{equation*}
Now, putting together the two previous inequalities yields;
\begin{equation}
-(\bw'(t),\bw(t))-\nu||\bw(t)||^2-a||\bw(t)||^2-b(|\bu|^\alpha
\bu-|\widetilde{\bu}|^\alpha \widetilde{\bu},\bu-\widetilde{\bu})\ge
0~~,~ \label{eq:4-3}
\end{equation}
where $\bw(t)=\bu(t)-\widetilde{\bu}$\,. From the monotonicity of
$T(\zeta)=|\zeta|^\alpha\zeta$, (\ref{eq:4-3}) imply that
\begin{equation*}
\displaystyle \frac{d}{dt}||\bw(t)||^2+2(\nu+a)||\bw(t)||^2\le 0~,
\end{equation*}
from which the announced estimate is readily obtained via
Gronwall's lemma. \QED

{\bf Acknowledgment} The authors are indebted to Prof Mamadou Sango
for interesting discussions on Brinkman equations. The work of
second author has been sponsored by the University of Pretoria and the National Research Foundation of
South Africa.

\bibliographystyle{plain}

\end{document}